\documentclass[12pt,reqno]{article}

\usepackage{cite,epic,eepic,euscript,verbatim,amsmath,amsthm,amssymb,amsfonts,afterpage,float,bm,stmaryrd,paralist,cancel,color,fancyhdr,graphics,graphicx,epsf,hyperref,makeidx,epsfig}

\usepackage{subfigure,epstopdf}
\usepackage{array}
\numberwithin{equation}{section}
\numberwithin{table}{section}
\numberwithin{figure}{section}

\newcommand{\reff}[1]{{\rm (\ref{#1})}}

\newcommand{\ve}{\varepsilon}          

\def\XXint#1#2#3{{\setbox0=\hbox{$#1{#2#3}{\int}$}
\vcenter{\hbox{$#2#3$}}\kern-.51\wd0}}

\newcommand{\td}{\tilde}

\begin{document}


\title
{
Computational Study on Hysteresis of Ion Channels: Multiple Solutions to Steady-State Poisson--Nernst--Planck Equations
}

\author{Jie Ding\thanks{
Department of Mathematics and Mathematical Center for Interdiscipline Research, Soochow University, 1 Shizi Street, Suzhou 215006, Jiangsu, China}
\and
Hui Sun\thanks{
Department of Mathematics and Statistics, California State University, Long Beach, CA, U. S. A.}
\and
Zhongming Wang\thanks{
Department of Mathematics and Statistics, Florida International University, Miami, FL, U. S. A.}
\and
Shenggao Zhou\thanks{
Department of Mathematics and Mathematical Center for Interdiscipline Research, Soochow University, 1 Shizi Street, Suzhou 215006, Jiangsu, China.
To whom correspondence should be addressed. E-mail: sgzhou@suda.edu.cn.
}
}

\date{}

\maketitle

\begin{abstract}
The steady-state Poisson--Nernst--Planck (ssPNP) equations are an effective model for the description of ionic transport in ion channels. It is observed that an ion channel exhibits voltage-dependent switching between open and closed states. Different conductance states of a channel imply that the ssPNP equations probably have multiple solutions with different level of currents. We propose numerical approaches to study multiple solutions to the ssPNP equations with multiple ionic species. To find complete current-voltage ($I$-$V$) and current-concentration ($I$-$C$) curves, we reformulate the ssPNP equations into four different  boundary value problems (BVPs). Numerical continuation approaches are developed to provide good initial guesses for iteratively solving algebraic equations resulting from discretization. Numerical continuations on $V$, $I$, and boundary concentrations result in S-shaped and double S-shaped ($I$-$V$ and $I$-$C$) curves for the ssPNP equations with multiple species of ions. There are five solutions to the ssPNP equations with five ionic species, when an applied voltage is given in certain intervals. Remarkably, the current through ion channels responds hysteretically to varying applied voltages and boundary concentrations, showing a memory effect.  In addition, we propose a useful computational approach to locate turning points of an $I$-$V$ curve.  With obtained locations, we are able to determine critical threshold values for hysteresis to occur and the interval for $V$ in which the ssPNP equations have multiple solutions. Our numerical results indicate that the developed numerical approaches have a promising potential in studying hysteretic  conductance states of ion channels.

 \bigskip

\end{abstract}

{ \allowdisplaybreaks
\section{Introduction}
\label{s:Introduction}
Essential for life, ion channels are protein molecules with a narrow spanning pore across membranes, regulating various crucial biological functions\cite{Hille_Book2001,MacKinnon04_ACIE04,IonChanel_HandbookCRC15}. They play fundamental roles in exchanging ions across cell membranes, propagating electric impulses in nerves, and maintaining
excitability of membrane\cite{NeuroBook,HH1_JPhys52}. Ion channels switch their conductance states  in response to variation of transmembrane voltages as well as ionic concentrations of cells and extracellular medium\cite{Sigworth_QRB94,Yamoah_BioPhyJ03,FologeaBBActa_2011, Bezrukov_EBioPhyJ15 ,Usherwood_Nat81,Cui_PflugersArch94,Nache_NatComm2013}.  For instance, lysenin channels that are inserted into a planar bilayer lipid membrane show voltage regulations with slowly changing external voltages \cite{FologeaBBActa_2011}. Voltage-dependent anion channels exhibit hysteretic response to a varying voltage with frequency of different magnitudes\cite{Bezrukov_EBioPhyJ15}. Ionic concentrations also have significant impacts on the switching of conductance states of voltage-gated channels. An L-type voltage-gated calcium channel shows distinct gating modes when different concentrations of charge carriers are applied\cite{Yamoah_BioPhyJ03}.

Hysteresis phenomenon is ubiquitous in optical devices\cite{Gibbs_Book}, many-body systems\cite{Hyst_RevModPhys99}, and biological systems\cite{Hyst_Math_Book}. In recent years, there has been a growing interest in understanding hysteretic response of ion channels to varying applied voltages\cite{Roope_JGPhys05,Bezrukov_JCP06,Andersson_MathBiosci2010,FologeaBBActa_2011,Das_PRE2012,Krueger_BiophyChem13, Bezrukov_EBioPhyJ15, Cuello_PNAS17}.  Hysteresis of ions channels is of physiological significance, since it involves many human physiological processes\cite{Roope_JGPhys05, Das_PRE2012, Cuello_PNAS17}. One distinct feature of hysteresis is the memory effect when the system undergoes transitions between different states.  In the context of voltage-gated ion channels, the current through channels increases and decreases along different paths when applied voltages periodically ascend and descend, respectively\cite{Bezrukov_JCP06,Andersson_MathBiosci2010,FologeaBBActa_2011,Das_PRE2012, Bezrukov_EBioPhyJ15}. It is pointed out that hysteresis takes place when frequency of an applied oscillating voltage is competing to typical relaxation time of transitions between different conductance states\cite{Bezrukov_JCP06,FologeaBBActa_2011,Das_PRE2012}. To explore such a hysteretic response, several discrete state Markov models have been developed\cite{HH1_JPhys52, HH2_JPhys52,HH3_JPhys52,Hodgkin_PRSL58,Altomare_JGenPhy01,Roope_JGPhys05,Bezrukov_JCP06,Noble_PTRSocA_09, Das_PRE2012}. In such models, ion channels are assumed to have certain number of states,  representing closed and open states. In addition, Markovian properties are assumed in the transitions between different states with certain transition rates. The master equation of stochastic processes is derived to describe the probability of the channel being in each state with respect to time.  It should be noted that hysteresis exhibited in voltage-gated ion channels is often associated with stochastic conformational changes of ion channel proteins. In this work, we also observe hysteretic response of currents to the varying applied voltages as well as ionic concentrations, using a deterministic model of ionic transport, rather than a stochastic description of different conductance states. Our results imply that, in addition to its stochastic nature, the gating phenomenon may possibly have deterministic factors associated to multiple states of ionic transport current through open channels.

The Poisson--Nernst--Planck (PNP) equations are an effective theory in modeling ionic transport through ion channels under electrostatic potential differences across a membrane. The Poisson equation describes electrostatic potential due to the charge density that stems both from mobile ions and fixed charges in the system. Nernst-Planck equations govern the diffusion and migration of ions in gradients of ionic concentrations and electrostatic potential. More complicated models have been developed to account for ionic steric effects and ion-ion correlations that are neglected in mean-field derivation of the PNP theory\cite{MaXu_JCP14,XuMaLiu_PRE14, ZhouWangLi_PRE11,LiLiuXuZhou_Nonliearity13, LiWenZhou_CMS16, BZLu_BiophyJ11, BZLu_JCP14, BZLu_JSP16, LinBob_CMS14,LeeHyonLinLiu_Nonlinearity11,HyonLiuBob_JPCB12,HyonLiuBob_CMS10,BobHyonLiu_JCP10,XuShengLiu_CMS14}. Due to nonlinear coupling of electrostatic potential and ionic concentrations, it is very hard to solve the problem analytically. Most of mathematical analyses of the PNP equations are based on singular perturbation methods, in which the Debye length is assumed to be much smaller than the dimensions of ion channels, giving rise to a small singular parameter. Such singular perturbation problems are solved mainly by two categories of methods: matched asymptotic expansions \cite{RubinsteinSIAM_Book,BarChenBob_SIAP92,BarChenBob_SIAP97,SchussNadlerBob_PRE01,SingerGillBob_ESIAM08,WangHeHuang_PRE14,BobLiu_SIADS08,BazantSteric_PRE07,BazantChuBayly_SIAP06} and geometric singular perturbation theory \cite{WLiu_SIAP05, BergLiu_SIMA07, WLiu_JDE09,LinLiuZhang_SIADS13,JiLiuZhang_SIAP15,BobLiuXu_Nonlinearity15,LiuXu_JDE15,JiaLIuZhang_DCDSB16}. For instance, Wang \emph{et al.} study the PNP equations using matched asymptotic analysis, and prove the existence and uniqueness of the solution to the PNP equations with two and three ionic species\cite{WangHeHuang_PRE14}. Using a geometric framework, Liu considers the steady-state PNP (ssPNP) equations with multiple ionic species, and reduces the problem to nonlinear algebraic equations\cite{WLiu_JDE09}. It is pointed out that the ssPNP equations with multiple ionic species pose a significantly more challenging task in studying their solutions analytically and even numerically.

The PNP equations are also known as the drift-diffusion model in the literature of semiconductor physics\cite{PMarkowich_Book}. Instead of ions that possibly have  more than two species, two types of particles, electrons and holes, are considered in semiconductor devices. Many global existence results on the solutions to steady-state drift-diffusion equation have been established\cite{PMarkowich_Book,RubinsteinSIAM_Book}, whereas results on the uniqueness of solutions are mainly limited to very small applied voltages.  The existence of multiple solutions to the one-dimensional steady-state drift-diffusion equation is mainly evidenced by numerical simulations. Mock designs an example, with piecewise constant fixed charges of large magnitude, that has multiple numerical solutions\cite{MockExample}. Multiple steady states to the drift-diffusion equation are also studied numerically and asymptotically, based on local electroneutrality approximations \cite{Rubinstein_SIAP87, Steinruck_SIAP89} However, it is found that asymptotic analysis based on local electroneutrality approximations is not so accurate within the parameter range for which the current-voltage relation is nonmonotonic\cite{Ward_SIAP91}.

The study of multiple steady states to generalized PNP-type models has also attracted much attention in recent years\cite{LinBob_CMS14,LinBob_Nonlinearity15,HungMihn_arXiv15,Gavish_arXiv17}. It is shown that,\cite{LinBob_Nonlinearity15,HungMihn_arXiv15} when zero-current boundary conditions are used, the steady-state Nernst-Planck-type equations can be further reduced to a system of algebraic equations, which defines generalized Boltzmann distributions, i.e.,functions of concentration against electrostatic potential. Rigorous analyses prove that there are multiple steady states to the PNP-type models with multiple ionic species and fixed charges. In a recent work\cite{Gavish_arXiv17}, Gavish considers a PNP-type model with two ionic species, and studies the solution trajectories in a phase plane composed of two ionic concentrations.  In this work, we numerically study multiple steady states of the classical PNP equations, with multiple ionic species, that have non-zero current through ion channels.

\begin{figure}[htbp]
    \centering
    \includegraphics[width=3.1in,height=2.5in]{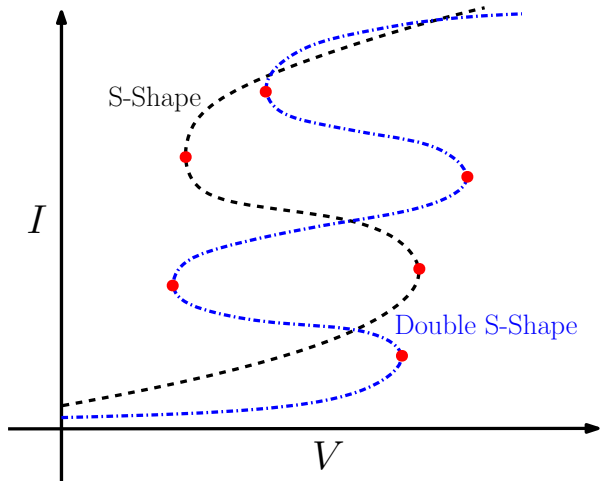}
    \caption{A schematic plot of $I$-$V$ characteristic curves. There are two types of curves: S-shape (in black) and double S-shape (in blue). The red dots represent turning points where $d V/d I =0$. } \label{f:Schem}
\end{figure}
The current-voltage ($I$-$V$) characteristic relation is often studied in the literature of ion channels. The total current ($I$) of an ion channel consists of contributions from each ionic species, and the voltage ($V$) is the potential difference between two ends of an channel.  The S-shaped $I$-$V$ curve shown in Figure\ref{f:Schem} has been found by solving steady-state PNP (ssPNP) equations in the literature of semiconductor physics\cite{MockExample, Rubinstein_SIAP87, Steinruck_SIAP89,Ward_SIAP91}.  Rather than two types of charged particles (electrons and holes), we here consider more than two ionic species in the context of ion channels. When five species of ions present (e.g., a mixture of $\mbox{K}^+$, $\mbox{Na}^+$, $\mbox{Ca}^{2+}$, $\mbox{Cl}^-$, and $\mbox{CO}_3^{2-}$), we find a double S-shaped $I$-$V$ curve for which there are five current levels of conductance when $V$ is given in some interval, cf. Figure~\ref{f:Schem}. To find these $I$-$V$ curves, we reformulate the ssPNP equations into four different boundary value problems (BVPs), which are solved with the help of newly developed numerical continuation approaches.  Another contribution of this work is the development of a computational method for locating turning points on a $I$-$V$ curve, i.e., the red dots in Figure~\ref{f:Schem}. Turning points are of great practical significance, since they not only give critical threshold values for hysteresis to occur but also specify the interval of applied voltages in which the multiple solutions present. In addition, we study the impact of boundary concentrations in cells and extracellular medium on currents through channels. Remarkably, we find that  currents respond hysteretically to variation of concentrations , resulting in an S-shaped curve in the current-concentration plane. To the best of our knowledge, such a hysteretic response of currents to concentrations has not been reported in the literature of PNP equations. Overall, our results imply that, within certain range of parameters, conductance states of ion channels depend on applied voltages and ionic concentrations as well as their history values, displaying a memory effect.

The rest of the paper is organized as follows. In section \ref{s:ModelDescrip}, we describe the ssPNP equations and reformulate them into four different boundary value problems.  In section \ref{s:ComputMethods}, we develop numerical continuation approaches and methods for computing  turning points on $I$-$V$ curves. Finally, the section \ref{s:NumRes} is devoted to showing our numerical results on multiple solutions of the ssPNP equations, turning points, and hysteretic response of currents to variation of applied voltages and concentrations.

\section{Model Description}\label{s:ModelDescrip}
\subsection{Governing Equations}
The transport of diffusive ionic or molecular species can be described by
\[
\partial_t c_i  + \nabla \cdot {\bf J}_i=0 \quad \mbox{for}~ i=1, \dots, M,
\]
where $c_1, \cdots, c_M$ and ${\bf J}_1, \cdots, {\bf J}_M$ are concentrations and flux densities of charged species, respectively. The flux of the $i$th species is defined by
\[
{\bf J}_i=-D_i \left(\nabla c_i + \beta z_i e c_i \nabla \psi \right),
\]
where $e$ is the elementary charge, $\psi$ is the electrostatic potential, $\beta$ is the inverse of thermal energy, and $D_i$ and $z_i$ are the diffusion constant and valence of $i$th species, respectively. The electrostatic potential is governed by the Poisson equation
\[
-\nabla \cdot \ve_0 \ve_r \nabla \psi = \rho,
\]
where $\ve_0$ is the vacuum permittivity and $\ve_r$ is the relative permittivity (or dielectric coefficient). The total charge density in the system consists of mobile ions and fixed charge:
\[
\rho=\sum_{i=1}^M z_i e c_i +\rho^f.
\]

Let $L$, $D_0$, and $c_0$ be the characteristic length, diffusion constant, and concentration, respectively. Introduce another  characteristic length $\lambda_D= \sqrt{\frac{\ve_0 \ve_r}{2 \beta e^2 c_0}}$ for an aqueous solution with bulk ionic concentration $c_0$ and homogenous dielectric coefficient $\ve_r$. We shall introduce the following dimensionless parameters and variables:
\begin{equation}\label{Rescale}
\td x = x/L, \td t =t D_0/L\lambda_D, \td c_i =c_i/c_0,  \td D_i =D_i/D_0, \td \rho^f= \rho^f/c_0, ~\mbox{and}~  \phi=\beta e \psi.
\end{equation}
Combining above equations and dropping all the tildes lead to nondimensionalized Poisson--Nernst--Planck (PNP) equations
\begin{equation}
\left\{
\begin{aligned}
&\partial_t c_i= \frac{\lambda_D}{L} \nabla \cdot D_i \left ( \nabla c_i + z_i c_i \nabla \phi \right),~ i=1, \dots, M,\\
&- 2 \frac{\lambda_D^2}{L^2}  \Delta \phi =\sum_{i=1} ^M z_i c_i + \rho^f.
\end{aligned}
\right.
\end{equation}

In this work, we consider steady states of a one-dimensional problem that describes the transport of ions through an ion channel with length $2L$. The computational domain is further rescaled to $[-1, 1]$, and the steady state PNP (ssPNP) equations are reduced to
\begin{equation}\label{PNP1D}
\left\{
\begin{aligned}
&\partial_x  \left ( \partial_x c_i + z_i c_i \partial_x \phi \right) = 0,~ i=1, \dots, M,\\
&- \partial_{xx} \phi =\kappa \left( \sum_{i=1} ^M z_i c_i + \rho^f \right),
\end{aligned}
\right.
\end{equation}
where $\kappa= \frac{L^2}{2\lambda_D^2}$ is a dimensionless parameter.

\subsection{Problem Formulations}\label{ss:ProbForms}
We are interested in the case that  two ends of an ion channel (or nanopore) are connected to individual ionic reservoirs, i.e.,
\begin{equation}\label{ConBcs}
c_i(-1)=c_i^L ~\mbox{and}~ c_i(1)=c_i^R ~\mbox{for} ~  i=1, \dots, M.
\end{equation}
Such concentrations satisfy the neutrality conditions
\begin{equation}\label{NeuCon}
\sum_{i=1}^M z_i c_i^L = \sum_{i=1}^M z_i c_i^R = 0.
\end{equation}
To describe electrostatic potential differences between two ends, we prescribe
\begin{equation}\label{PotBCs}
\phi(-1)=0 ~\mbox{and}~ \phi(1)=V,
\end{equation}
where $V$ is the applied potential difference. We define the total current
\[
I = \sum_{i=1}^M z_i J_i,
\]
where $J_i=\partial_x c_i + z_i c_i \partial_x \phi$.  With a given applied potential difference $V$, we solve the ssPNP equations \reff{PNP1D} with boundary conditions \reff{ConBcs} and \reff{PotBCs}, and obtain the total current $I$ across the channel. By varying applied potential differences, we get the current-voltage relation ($I$-$V$ curve). Let $\mu=\phi'$, ${\bf Y} =(\phi, \mu, c_1, J_1, \cdots, c_M, J_M)^T$, and ${\bf F}= ( \mu, -\kappa\left(  \sum_{i=1}^M z_i c_i + \rho^f\right), J_1-z_1 c_1 \mu, 0, J_2-z_2 c_2 \mu, 0, \cdots, J_M-z_M c_M \mu, 0)^T$. Our problem can be written as
\[
\mbox{Problem (V2I)}:  ~~{\bf Y}' = {\bf F} ({\bf Y})\qquad \mbox{and} \qquad
\left\{
\begin{aligned}
&\phi (-1)= 0,\\
&\phi (1)= V,\\
&c_i(-1)= c_i^L, ~i=1, \dots, M,\\
&c_i(1)= c_i^R, ~~i=1, \dots, M.\\
\end{aligned}
\right.
\]
Notice that we call this formulation ``Problem (V2I)'' throughout the following contents.

As shown in Figure~\ref{f:Schem}, the solution to Problem (V2I) may not be unique for some applied voltages. There are three solutions for the S-shaped $I$-$V$ curve and five solutions for the double S-shaped $I$-$V$ curve, when $V$ is prescribed in some interval. The Problem (V2I) could be a large nonlinear system if multiple species of ions are considered. When numerically solving such a  system, good initial guesses are crucial to devising convergent iterations. Numerical continuation is a powerful tool of providing initial guesses. See section \ref{ss:Cont} for more details. Continuation on applied voltages helps in finding the the low-current branch and high-current branch for S-shaped and double S-shaped curves. However, continuation on $V$ often misses finding the intermediate branches when multiple solutions present, c.f.,~Figure\ref{f:2ionsIV} and Figure\ref{f:5ionsVCont}. In addition,  the continuation advances with very small stepsizes as $V$ approaches the turning points where the Jacobian of the discretized nonlinear system becomes more and more singular.  In order to overcome these drawbacks, we view the current $V$ as a function of $I$, which is single valued. To find a complete characteristic curve in the $I$-$V$ plane, we prescribe $I$ in the computation, instead of $V$, and solve the following problem with continuations on $I$:
 \[
\mbox{Problem (I2V)}:  ~~{\bf Y}' = {\bf F} ({\bf Y})\qquad \mbox{and} \qquad
\left\{
\begin{aligned}
&\phi (-1)= 0,\\
&z_1 J_1 (1) + \cdots +z_M J_M(1) =I,\\
&c_i(-1)= c_i^L, ~i=1, \dots, M,\\
&c_i(1)= c_i^R, ~~i=1, \dots, M.\\
\end{aligned}
\right.
\]
Similarly, we notice that this formulation is called  ``Problem (I2V)'' throughout the following contents.  We obtain the $I$-$V$ curve by collecting the voltage at the right boundary, i.e., $V=\phi(1)$.

We are also interested in understanding the effect of boundary concentrations \reff{ConBcs} on the currents through a channel. For simplicity, we fix concentration of each species of at the left end and vary the concentration of two species of ions with indice $1$ and $2$ at the right end. By relabeling the ions, we can assume that these two species of ions have opposite signs, to satisfy the neutrality conditions \reff{NeuCon}. We define the following problem
\[
\mbox{Problem (C2I)}:  ~~{\bf Y}' = {\bf F} ({\bf Y})\qquad \mbox{and} \qquad
\left\{
\begin{aligned}
&\phi (-1)= 0,\\
&\phi (1)= V,\\
&c_i(-1)= c_i^L, ~i=1, \dots, M,\\
&c_1(1)= c_B,\\
&c_2(1)= -\left( z_1c_B  + z_3  c_3^R + \cdots +z_M c_M^R  \right)/z_2, \\
&c_i(1)= c_i^R, ~~i=3, \dots, M,\\
\end{aligned}
\right.
\]
where we introduce a variable $c_B$. We solve this problem with varying $c_B$ and compute the corresponding current $I$. Similar to the case of $I$-$V$ relation, it is helpful to define the following problem
 \[
\mbox{Problem (I2C)}:  ~~{\bf Y}' = {\bf F} ({\bf Y})\qquad \mbox{and} \qquad
\left\{
\begin{aligned}
&\phi (-1)= 0,\\
&\phi (1)= V,\\
&c_i(-1)= c_i^L, ~i=1, \dots, M,\\
&z_1 J_1 (1) + \cdots +z_M J_M(1) =I,\\
&z_1c_1(1)+z_2c_2(1)= - z_3 c_3^R  - \cdots -z_M c_M^R  , \\
&c_i(1)= c_i^R, ~~i=3, \dots, M,\\
\end{aligned}
\right.
\]
where we include the neutrality conditions \reff{NeuCon} as a boundary condition for concentrations. We solve this problem with varying $I$ and compute $c_1$ at right boundary to get $c_B$, i.e., $c_B=c_1(1)$. See numerical results shown in Figure~\ref{f:2ionsIC} and corresponding descriptions for more details.
\section{Computational Methods}\label{s:ComputMethods}
 The boundary value problems (BVPs) defined above are numerically solved with the program of BVP4C \cite{BVP4C} in the Matlab. One of its advantages is that the algorithm uses a collocation method to discretize a BVP on an adaptive, nonuniform mesh. Grid points are more densely distributed at locations where the solutions have large variations, cf.~Fig \ref{f:GPD}. The resulting nonlinear algebraic equations are solved iteratively by linearizations, e.g., Newton-type methods. Often analytical Jacobians  $\partial {\bf F}/\partial {\bf Y}$ are provided to accelerate the iterations. For Problems  (V2I), (I2V), (C2I), and (I2C), we have the following Jacobian
\[
\frac{\partial {\bf F}}{\partial {\bf Y}}=
\left(
\begin{array}{ccccccccc}
0 & 1 & 0               & 0 & 0                  &0 &\cdots &0                 &0 \\
0 & 0 & -\kappa z_1& 0 & -\kappa z_2 &0 &\cdots &-\kappa z_M &0 \\
0 & -z_1c_1 & -z_1\mu& 1 & 0 & 0 &\cdots &0 & 0 \\
0 & 0  & 0                & 0 & 0                  &0 &\cdots &0                &0 \\
0 & -z_2c_2 & 0 & 0 & -z_2\mu& 1 &\cdots & 0 & 0 \\
0 & 0  & 0                & 0 & 0                  &0 &\cdots &0                &0 \\
   &   & \cdots & &  & &\cdots & &\\
      &   & \cdots & &  & &\cdots & &\\
0 & -z_Mc_M & 0 & 0 &0 & 0 &\cdots &  -z_M\mu & 1 \\
0 & 0  & 0                & 0 & 0                  &0 &\cdots &0                &0 \\
\end{array}
\right).
\]
In our computations, we take the following piecewise constant profile for fixed charges\cite{MockExample, Rubinstein_SIAP87, Steinruck_SIAP89, Ward_SIAP91}:
\begin{equation}\label{rhs}
\rho^f(x)= \sigma \rho_i^f ~\mbox{for}~ x \in \left (x_{i-1}, x_i \right),
\end{equation}
where $i=1, \dots, N$, $\sigma >0$, $x_0=-1$, $x_i=x_{i-1}+L_i$, and $\sum_{i=1}^N L_i = 2$. Such a particular profile of fixed charges has its practical meaning in the context of ion channels. Along a channel, fixed charges carried by atoms in the membrane could carry positive or negative partial charges. In our current treatment, the values of these partial charges are simply approximated by constants of alternating signs. 
\subsection{Strategy of Continuation}\label{ss:Cont}
As mentioned above, BVP4C solves nonlinear algebraic equations, resulting from discretization, by Newton-type iterative methods. It is well known that the convergence of Newton's iterations highly depends on the the choice of initial guesses. In our problems, the main difficulty arises from the fixed charges $\rho^f$ that changes its sign drastically. Therefore, we adopt a strategy of continuation on the parameter $\sigma$, cf. \reff{rhs}. To be specific, we gradually increase the value of $\sigma$ with nonuniform increments, i.e.,
\[
0=\sigma_0 < \sigma_1 <\sigma_2<\cdots < \sigma_K = \sigma.
\]
For the $k$th stage of iteration, we solve the problem with $\sigma_k$ using the solution from the $(k-1)$th stage with $\sigma_{k-1}$ as an initial guess.  In each stage, the initial guess is close to final solutions and therefore  iterations converge easily.

To obtain $I$-$V$ curves and current-concentration relations, we also apply continuation approaches on $V$, $I$, and $c_B$ when the problems (V2I), (I2V), (C2I), and (I2C) are considered. When using continuation on $V$ for Problem (V2I), care should be taken as $V$ approaches the turning points at which the convergence of Newton's iterations significantly slows down and Jacobians of iterations become more and more singular. We decrease the increments of continuation appropriately at turning points. Similar observations can be made for the continuation on $c_B$. See Figure~\ref{f:2ionsIV}, \ref{f:2ionsIC}, \ref{f:5ionsVCont}, and related descriptions for more details.
\subsection{Turning Points on $I$-$V$ Curves}\label{ss:TPMethods}
It is of practical importance to study the location of turning points, due to the fact that the $V$ values of these critical points are threshold values for hysteresis to take place, and that they are endpoints of intervals for $V$ in which the ssPNP equations have multiple solutions.  Turning points are also locations where continuation increments on $V$ should be carefully chosen,  when solving Problem (V2I) with Newton's iterations. It is observed from Figure~\ref{f:Schem} that, if $V$ is viewed as a differentiable function of $I$, the equation $d V/d I =0$ holds at turning points (labeled with red dots) of an $I$-$V$ curve.

To incorporate the information of $d V/d I$, we perform sensitivity analysis on the ssPNP equations with respect to variation of $I$. We first introduce the current $I$ as an unknown and reformulate Problem (V2I) as follows:
\begin{equation}\label{ODEIForm}
\left\{
\begin{aligned}
&\phi' =\mu,\\
&\mu' =- \kappa \left(\sum_{j=1}^{M} z_j c_j + \rho^f \right),\\
&c_i'=J_i-z_i c_i \mu,~~ i=1, \dots, M-1,\\
&J_i'=0, ~~i=1, \dots, M-1,\\
&c_M'= \frac{I-\sum_{j=1}^{M-1} z_j J_j}{z_M}-z_M c_M \mu,\\
&I'=0,
\end{aligned}
\right.
\qquad \mbox{and} \qquad
\left\{
\begin{aligned}
&\phi (-1)= 0,\\
&\phi (1)= V,\\
&c_i(-1)= c_i^L, ~i=1, \dots, M,\\
&c_i(1)= c_i^R, ~~i=1, \dots, M,\\
\end{aligned}
\right.
\end{equation}
where we replace the unknown $J_M$ in Problem (V2I) and (I2V) by $I$. We introduce an infinitesimal perturbation to $I$, i.e., $I^\delta = I +\delta I$. Such a perturbation gives rise to perturbations in other quantities: $\phi^\delta = \phi+\delta \phi$,  $V^\delta =V+\delta V$, $\mu^\delta = \mu+\delta \mu$, $c_i^\delta = c_i+\delta c_i$, and $J_i^\delta= J_i+\delta J_i$ for $i=1, \dots, M$. Plugging  these perturbed quantities into \reff{ODEIForm} and subtracting from \reff{ODEIForm}, we have by neglecting high order terms that
\begin{equation}\label{PertODEIForm}
\left\{
\begin{aligned}
& \delta \phi'=\delta \mu,\\
&\delta \mu'=- \kappa \sum_{j=1}^{M}  z_j \delta c_j,\\
&\delta c_i'=\delta J_i-z_i \mu\delta c_i  -z_i c_i  \delta\mu,~~ i=1, \dots, M-1,\\
&\delta J_i'=0,~ i=1, \dots, M-1,\\
&\delta c_M'= \frac{\delta I-\sum_{j=1}^{M-1} z_j \delta J_j}{z_M}-z_M \mu \delta c_M - z_M c_M  \delta\mu ,\\
& \delta I'=0,
\end{aligned}
\right.
~ \mbox{and} ~
\left\{
\begin{aligned}
&\delta \phi (-1)= 0,\\
&\delta \phi (1)= \delta V,\\
&\delta c_i(-1)= 0, ~i=1, \dots, M,\\
&\delta c_i(1)= 0, ~~i=1, \dots, M.\\
\end{aligned}
\right.
\end{equation}
Dividing each equation in \reff{PertODEIForm} by $\delta I$, we have
\[
\left\{
\begin{aligned}
& \left(\frac{\delta\phi}{\delta I}\right)'=\frac{\delta \mu}{\delta I},\\
&\left(\frac{\delta\mu}{\delta I}\right)' =- \kappa \sum_{j=1}^{M}  z_j \frac{\delta c_j}{\delta I} ,\\
&\left(\frac{\delta c_i}{\delta I}\right)' =\frac{\delta J_i}{\delta I}-z_i \mu\frac{\delta c_i}{\delta I}  -z_i c_i  \frac{\delta \mu}{\delta I},~ i=1, \dots, M-1,\\
&\left(\frac{\delta J_i}{\delta I}\right)'=0,~ i=1, \dots, M-1,\\
&\left(\frac{\delta c_M}{\delta I}\right)' = \frac{1 -\sum_{j=1}^{M-1} z_j \frac{\delta J_i}{\delta I}}{z_M}-z_M \mu\frac{{\delta c}_M}{\delta I}  -z_M c_M  \frac{\delta \mu}{\delta I} ,
\end{aligned}
\right.
 \mbox{and} ~
\left\{
\begin{aligned}
&\frac{\delta\phi}{\delta I} (-1)= 0,\\
&\frac{\delta\phi}{\delta I} (1)= \frac{\delta V}{\delta I},\\
&\frac{\delta c_i}{\delta I}(-1)= 0, ~i=1, \dots, M,\\
&\frac{\delta c_i}{\delta I}(1)= 0, ~~i=1, \dots, M,\\
\end{aligned}
\right.
\]
where $\delta I$ is a constant with respect to $x$, by the fact that $\delta I' = \sum_{i=1}^{M} z_i J_i' =0 $. We denote by $\hat{\phi}=\displaystyle \lim_{\delta I \to 0} \delta\phi/\delta I$, $dV/d I=\displaystyle \lim_{\delta I \to 0} \delta V/\delta I$, $\hat{\mu}=\displaystyle \lim_{\delta I \to 0} \delta\mu /\delta I$, $ \hat{c}_i = \displaystyle \lim_{\delta I \to 0} \delta c_i/\delta I$, and $\hat{J}_i = \displaystyle \lim_{\delta I \to 0} \delta J_i/\delta I$ for $i=1, \dots, M$. As $\delta I$ goes to zero, we have by coupling \reff{ODEIForm} that
\[
\mbox{Problem (IVAug)}:
\left\{
\begin{aligned}
&\phi' =\mu,\\
&\mu' =- \kappa \left(\sum_{j=1}^{M} z_j c_j + \rho^f \right),\\
&c_i'=J_i-z_i c_i \mu,~~ i=1, \dots, M-1,\\
&J_i'=0, ~~i=1, \dots, M-1,\\
&c_M'= \frac{I-\sum_{j=1}^{M-1} z_j J_j}{z_M}-z_M c_M \mu,\\
&I'=0, \\
& \hat{\phi}'=\hat{\mu},\\
&\hat{\mu}' =- \kappa \sum_{j=1}^{M}  z_j \hat{c}_j ,\\
&\hat{c}_i' = \hat{J}_i-z_i \mu \hat{c}_i  -z_i c_i  \hat{\mu}, ~ i=1, \dots, M-1,\\
&\hat{J}_i'=0,~ i=1, \dots, M-1,\\
&\hat{c}_M' = \frac{1 -\sum_{j=1}^{M-1} z_j \hat{J}_j }{z_M}-z_M \mu\hat{c}_M  -z_M c_M  \hat{\mu},
\end{aligned}
\right.
 \hspace{-5mm} \mbox{and} ~
\left\{
\begin{aligned}
&\phi (-1)= 0,\\
&c_i(-1)= c_i^L, ~i=1, \dots, M,\\
&c_i(1)= c_i^R, ~~i=1, \dots, M,\\
&\hat{\phi} (-1)= 0,\\
&\hat{\phi} (1)= \frac{dV}{d I},\\
&\hat{c}_i(-1)= 0, ~i=1, \dots, M,\\
&\hat{c}_i(1)= 0, ~~i=1, \dots, M.\\
\end{aligned}
\right.
\]
Notice that this formulation is called  ``Problem (IVAug)'' throughout the following contents. We remark that there are $4M+3$ differential equations for $4M+3$ unknowns with $4M+3$ boundary conditions. Therefore, this problem determines a point on the $I$-$V$ curve with a given value of $d V/d I$.  As seen from Figure\ref{f:Schem}, this problem may have multiple solutions for an S-shaped or double S-shaped $I$-$V$ curve. Different initial guesses should be deliberately designed to achieve multiple solutions. It is of primary interest to calculate the the turning points of the $I$-$V$ curve where  $d V/d I=0$, cf. Figure~\ref{f:2ionsTurningPts}. After finding the turning points, we are able to determine the intervals for $V$ in which the ssPNP equations have multiple solutions.

\section{Numerical Results}\label{s:NumRes}
\begin{figure}[htbp]
    \centering
    \includegraphics[scale=0.62]{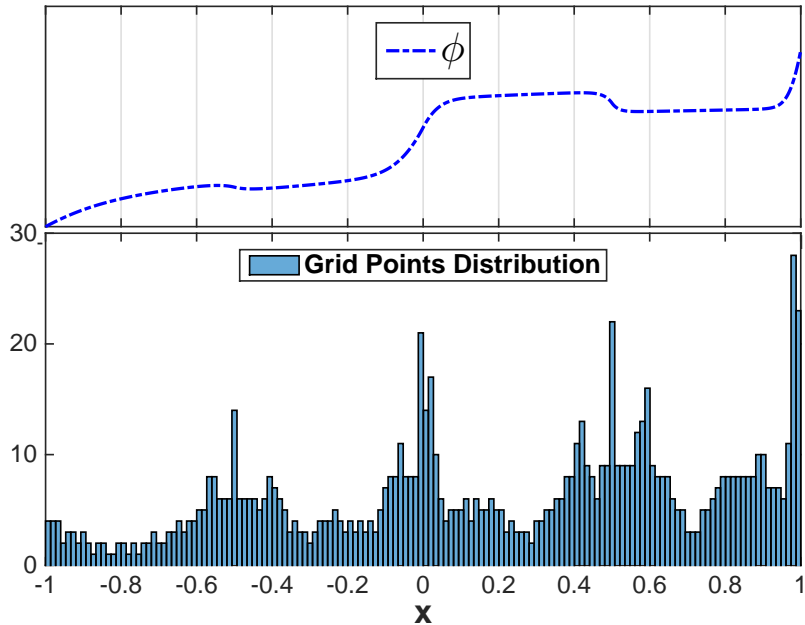}
    \caption{A typical solution of electrostatic potential and corresponding computational grid points distribution.} \label{f:GPD}
\end{figure}
In our computations with BVP4C, iterations for the discretized algebraic system stop when an absolute error is less than $10^{-6}$. The initial guesses for the very first calculation are provided by using a continuation approach on $\sigma$. We numerically solve the problem (V2I) using an adaptive, non-uniform mesh.   In Figure~\ref{f:GPD}, we show a typical solution of electrostatic potential and distribution of corresponding computational grid points. We can observe that  there are several peaks in the distribution when the solution has large variations at jumps of fixed charges, showing that the computational mesh  refines adaptively.

Unless otherwise stated, we use the following parameters in our computations: $\ve_r= 80$, $T=300$ K, $c_0=0.2$ M, and $\lambda_D=0.687$ nm. All the quantities that we show below are rescaled according to \reff{Rescale}. The values of the parameters, such as $\sigma$, $\kappa$, and $\rho_i^f$, are under normal physiological conditions of typical ion channels.

\subsection{Voltage-induced Hysteresis}\label{s:VHyst}
In this example, we study the response of current through an ion channel with two ionic species to varying applied voltages. We take the following values of parameters: $\kappa=60$, $\sigma=1$, $M=2$, $z_1=1$, $z_2=-1$, $c_1^L=1$,  $c_2^L=1$, $c_1^R=0.5$, $c_2^R=0.5$, $N=4$, $\rho_1^f=1$, $\rho_2^f =-10$, $\rho_3^f=20$, $\rho_4^f=-60$, and $L_i=0.5$ for $i=1, \dots, 4$.
\begin{figure}[htbp]
\centering
\includegraphics[scale=.7]{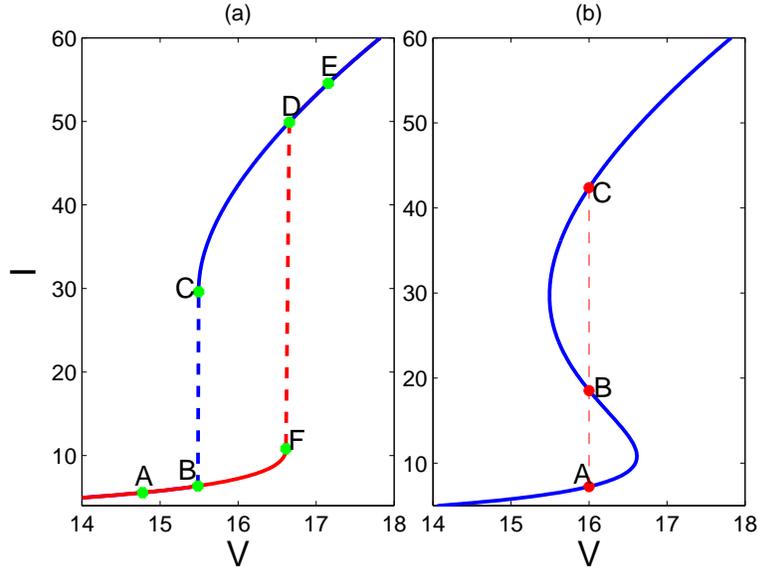}
\caption{ (a): Two branches of an $I$-$V$ curve obtained by solving the problem (V2I) with increasing and decreasing applied voltages; (b): A complete $I$-$V$ curve obtained by solving the problem (I2V) with increasing and decreasing $I$. When $V=16$, there are three $I$ values at red dots, labeled with A, B, and C. See Figure~\ref{f:2ions3Soluns} more details of each solution.}
 \label{f:2ionsIV}
\end{figure}

\begin{figure}[htbp]
\centering
\includegraphics[scale=.6]{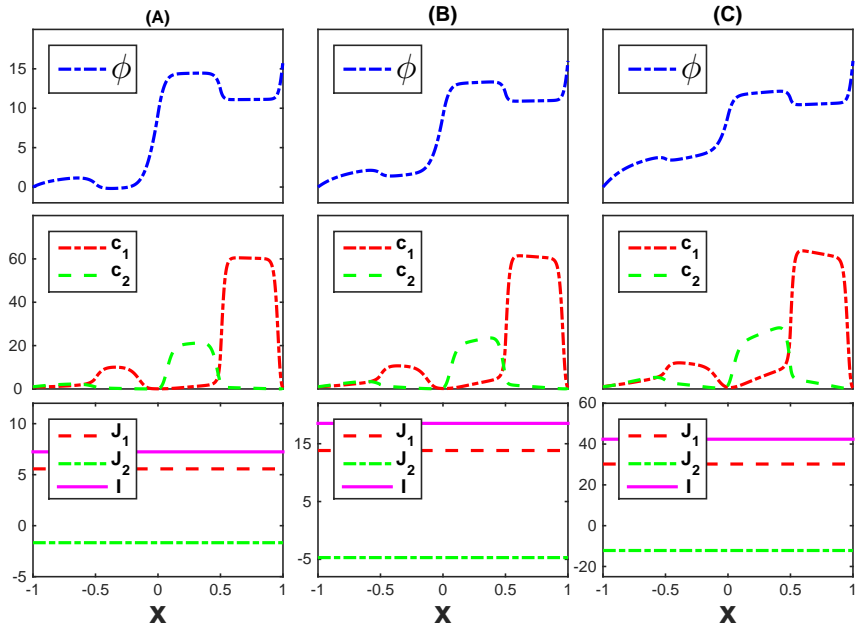}
\caption{ The solution profiles of electrostatic potential, concentrations, and currents, corresponding to A, B, and C in Figure~\ref{f:2ionsIV}(b). }
 \label{f:2ions3Soluns}
\end{figure}

We study a current-voltage relation with the help of numerical continuations proposed in section \ref{ss:Cont}. We first gradually increase an applied voltage and see the current ascending along a low-current branch, as shown in red in Figure~\ref{f:2ionsIV}(a). It is remarked that the convergence of iterations for solving the problem (V2I) gets more and more slow, as the applied voltage approaches a critical value ($V$ value of point F). When the applied voltage exceeds this critical value, the current all of a sudden jumps to a high-current (blue) branch, switching its conductance state. In a second round of sweeping, in which the voltage is decreased from large $V$ values, we observe that the current follows a totally different path ( from E to C) and suddenly jumps at C to a low-current branch, showing a typical hysteresis loop. We can see that the current of an ion channel depends on the applied voltage as well as its history values, exhibiting a memory effect. Although our model is purely deterministic, the sudden switching of conductance states is reminiscent of gating phenomenon. To explore more of the intermediate region of two branches, we view the voltage as a function of the current, and solve the problem (I2V) with continuations on $I$. It is noted that the iterations for solving this problem do not have any singular Jacobian as the curve passes points $C$ or $F$ in Figure~\ref{f:2ionsIV}(a).    Again, we design two rounds of sweeping with increasing and decreasing $I$, but the resulting plot of $I$-$V$ curves are almost identical, cf. Figure~\ref{f:2ionsIV}(b). Interestingly, there is an intermediate branch between the low-current and high-current branches. As such, there are three conductance states when $V$ belongs to some interval. For instance, there are three solutions to the problem (V2I) at red dots labeled with A, B, and C, when $V=16$. Figure~\ref{f:2ions3Soluns} displays profiles of electrostatic potential, concentrations, and currents, corresponding to these three points. We observe that the solutions of electrostatic potential resemble each other with the same potential differences, $V=16$. The concentrations for each case are all positive, meaning that three solutions are all of physical interests. As expected, the profiles of currents are constant and differ a lot for different conductance states. The current for each ionic species has larger magnitude when the channel is in a high conductance state.

\begin{figure}[htbp]
\centering
\subfigure{\includegraphics[scale=.55]{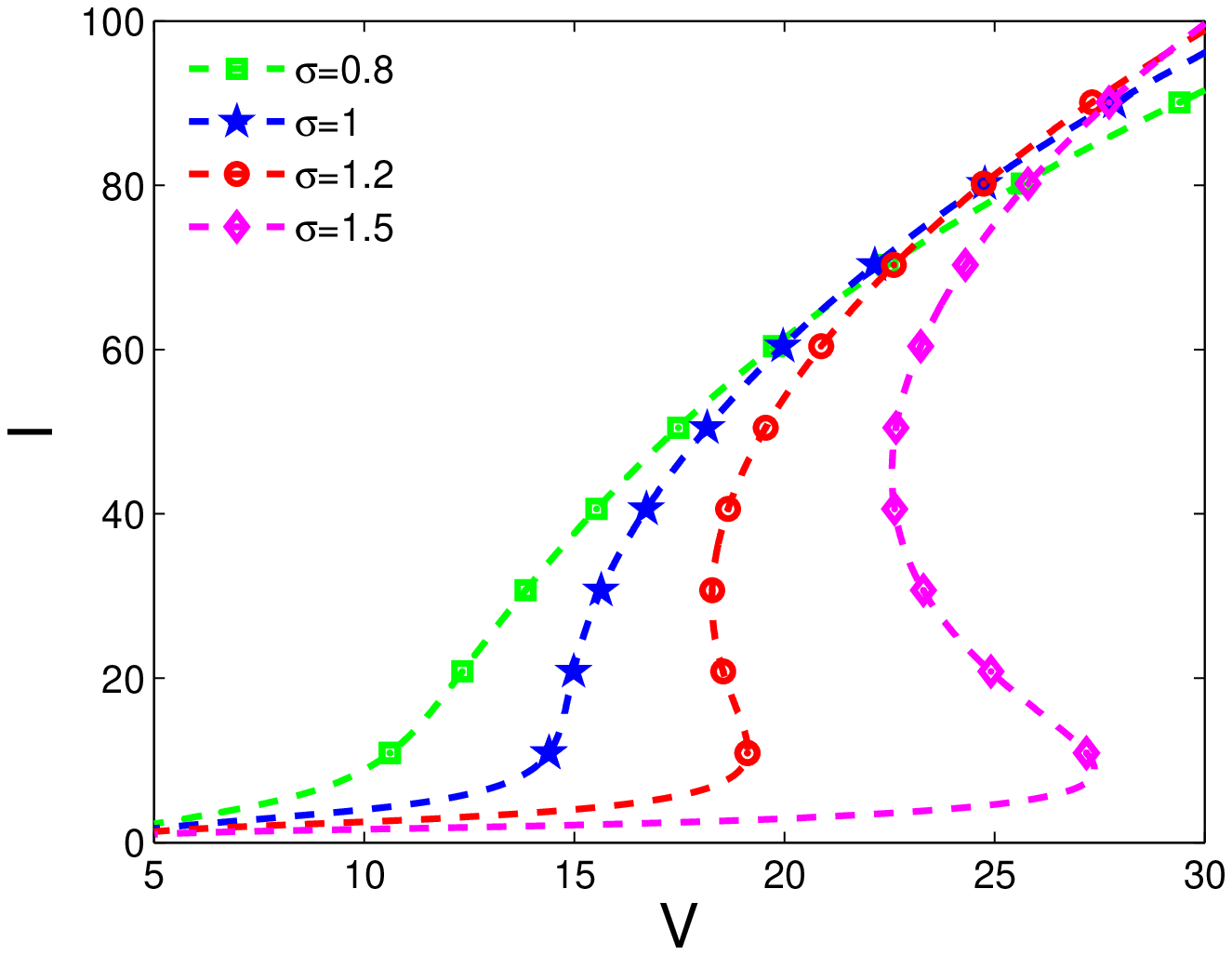}}
\subfigure{\includegraphics[scale=.55]{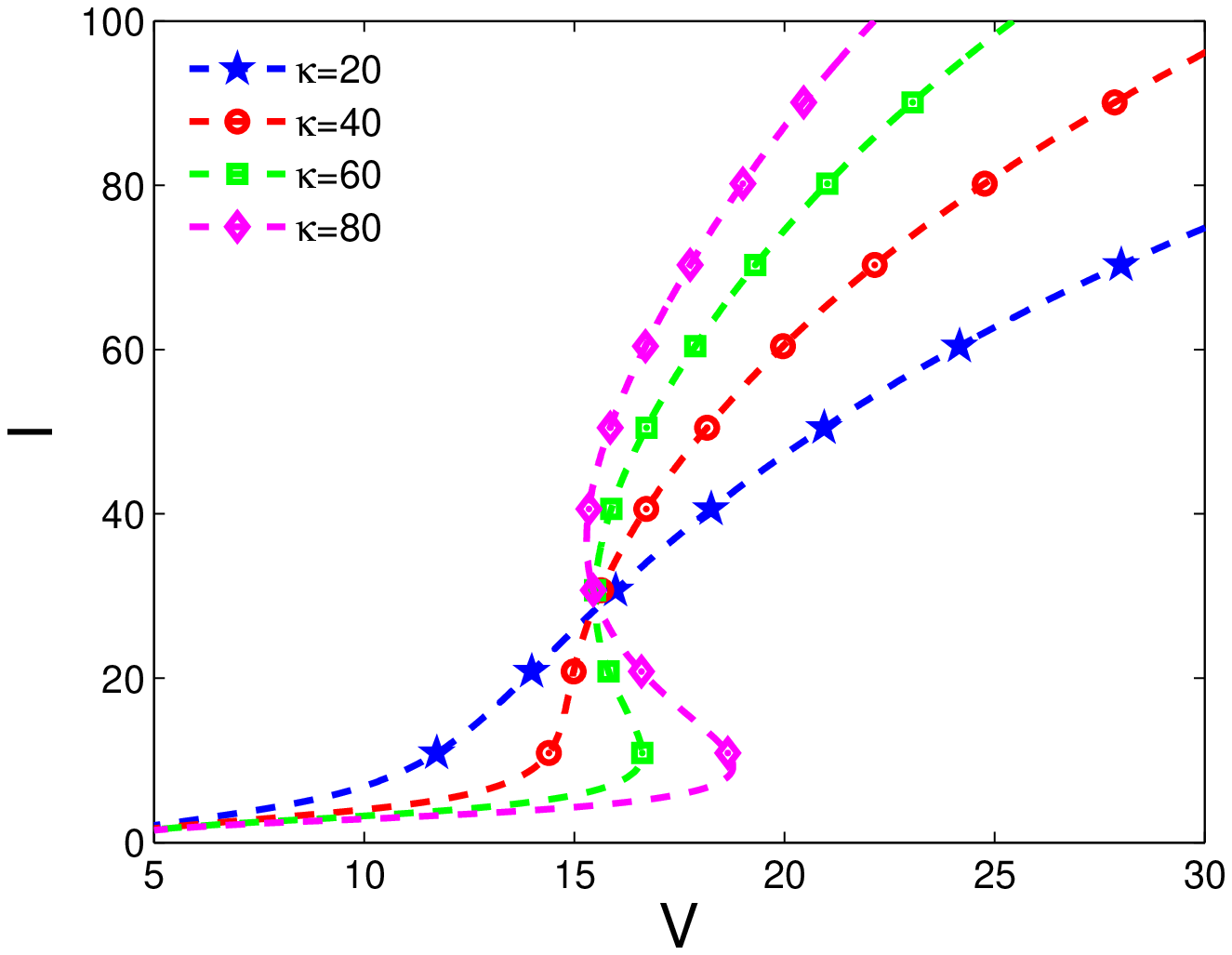}}
\caption{Left: $I$-$V$ curves for different $\sigma$ values when $\kappa=40$; Right:  $I$-$V$ curves for different $\kappa$ values when $\sigma =1$.}
 \label{f:2ionKappaRho4Fixed}
\end{figure}

To further understand hysteresis, we study the effect of parameters, $\kappa$ and $\sigma$, on the presence of multiple solutions. We see from Figure~\ref{f:2ionKappaRho4Fixed} that the values of $\kappa$ and $\sigma$ have significant impact on the shape of $I$-$V$ curve. We recall that $\sigma$, defined in \reff{rhs}, represents the magnitude of the fixed charge, and that $\kappa= \frac{L^2}{2\lambda_D^2}$ corresponds to the length of an ion channel under consideration.  We first set $\kappa=40$ and find that the curve gradually turns into a more and more obvious S shape, as $\sigma$ increases from $0.8$ to $1.5$. As such, an ion channel with larger magnitude of fixed charges is more likely to have multiple conductance states.   Also, we consider different values of $\kappa$ with fixed $\sigma=1$. It is interesting to see that the curve starts to switch back as the $\kappa$ exceeds $45$ or so, and that the $V$ value for which multiple solutions occur does not grow much as $\kappa$ increases. This manifests that, with other parameters fixed, it is easier to have multiple conductance states for a longer ion channel.

\begin{figure}[htbp]
\centering
\includegraphics[scale=.48]{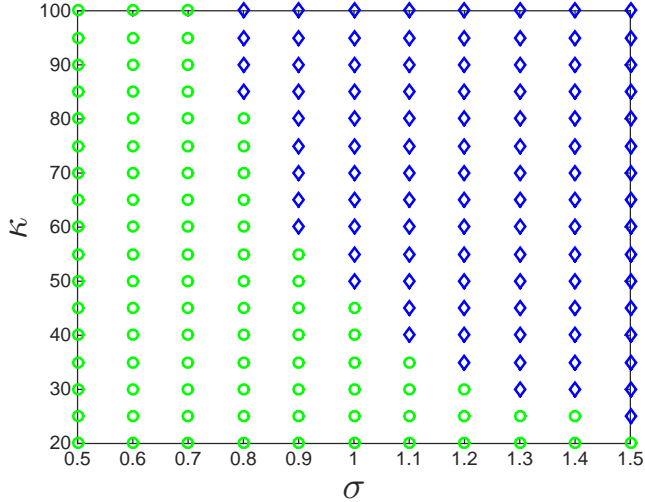}
\caption{A $\sigma$-$\kappa$ phase diagram showing monotonic and S-shaped $I$-$V$ curves. Green circles represent monotonic $I$-$V$ curves, and  blue diamonds represent S-shaped $I$-$V$ curves.}
\label{f:2ionsPhaseDiagram}
\end{figure}
As shown in Figure \ref{f:2ionsPhaseDiagram}, different combinations of $\sigma$ and $\kappa$ give rise to different shapes of $I$-$V$ curves. The transition value of $\kappa$ from a monotonic curve to a S-shaped one gets smaller as $\sigma$ increases. For larger $\sigma$, i.e., ion channels with larger magnitude of fix charges, it is more likely to have hysteretic responses to applied voltages.

\subsection{Concentration-induced Hysteresis}\label{s:CHyst}
In this example, we investigate the effect of boundary concentrations on the existence of multiple solutions to the ssPNP equations. For simplicity, we probe the solution behavior with a single varying variable, $c_B$, which is defined in section \ref{ss:ProbForms}. We take $\kappa=60$, $V=16$, and the same fixed charge profile as in the previous section.
\begin{figure}[htbp]
\centering
\includegraphics[scale=.7]{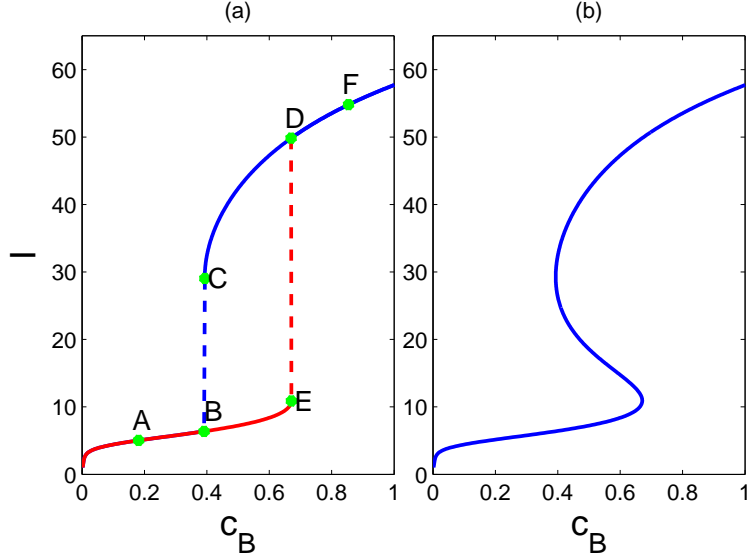}
\caption{ Current-concentration curves showing concentration-induced hysteresis.}
 \label{f:2ionsIC}
\end{figure}

We first solve the problem (C2I) using continuations with an increasing $c_B$ and find a low-current branch (red curve) in Figure~\ref{f:2ionsIC}. When $c_B$ goes over a critical value (concentration of $E$), the current jumps up to $D$ and follows a  high-current branch, shown in blue in Figure~\ref{f:2ionsIC} (a). Another round of sweeping starting from large $c_B$ values shows that the current stays in the high-current branch until $c_B$ gets less than the x-axis value of point C. This is very similar to the hysteresis loop described in previous section for $I$-$V$ curves. To find a complete current-concentration curve, we also conduct continuations on $I$ by solving the problem (I2C).  The curve has an intermediate branch that connects the points C and E. As such, we find three solutions, for some boundary concentrations, to the ssPNP equations with a fixed applied voltage. It is well documented that ionic concentrations have pronounced influence on the switching of conductance states of ion channels\cite{Usherwood_Nat81,Cui_PflugersArch94,Yamoah_BioPhyJ03, Nache_NatComm2013}.  To the best of our knowledge, such a hysteretic response of current to boundary concentrations has not been studied yet in the literature of PNP equations. In summary, we conclude that the existence of multiple solutions to the ssPNP equations depends sensitively on the coefficient $\kappa$, the fixed charges $\rho^f$, applied voltages $V$, and boundary concentrations.

\subsection{Double S-shaped $I$-$V$ curve}\label{s:SSIV}
\begin{figure}[htbp]
\centering
\includegraphics[scale=.75]{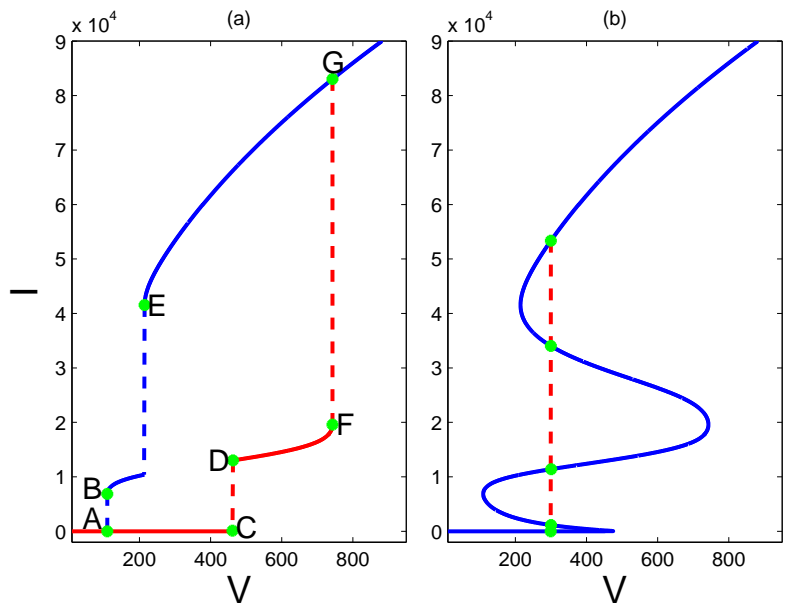}
\caption{(a): Several branches of an $I$-$V$ curve obtained by solving the problem (V2I) with increasing and decreasing applied voltages; (b): A complete $I$-$V$ curve obtained by solving the problem (I2V) with increasing and decreasing $I$. There are five solutions indicated by green dots. }
\label{f:5ionsVCont}
\end{figure}
In this example, we study conductance states of ion channels with five species of ions. As far as we know, the ssPNP equations with more than three ions have not been well studied numerically. In the literature of semiconductor physics, the PNP (drift-diffusion) equations are often studied with only two species, i.e., electrons and holes. It is noted that the computational complexity increases significantly as the number of species increases\cite{WLiu_JDE09, XuMaLiu_PRE14,LinBob_Nonlinearity15}.   We here take the following values of parameters:  $\kappa=200$, $\sigma=1$, $M=5$, $z_1=1$, $z_2=-1$, $z_3=2$, $z_4=-2$, $z_5=1$, $c_1^L=1$,  $c_2^L=1$, $c_3^L=0.5$,  $c_4^L=1$, $c_5^L=1$, $c_1^R=0.5$,  $c_2^R=2$, $c_3^R=1$,  $c_4^R=0.5$, $c_5^R=0.5$, $N=4$, $\rho_1^f=720$, $\rho_2^f =-800$, $\rho_3^f=960$, $\rho_4^f=-5600$, and $L_1=0.4$, $L_2=0.6$, $L_3=0.8$, $L_4=0.2$.

From Figure~\ref{f:5ionsVCont} (a), we can see that, as $V$ increases, the current increases gradually along the blue branch and suddenly jumps at $C$ and $F$ to branches of higher conductance. In a second round of sweeping with decreasing $V$,  we find that the current moves along a different branch and has abrupt jumps at $E$ and $B$, switching its conductance states. Similar to hysteresis loops, the conductance state of the channel has a memory effect, i.e., both the applied voltages and their history values influence the conductance state. As shown in  Figure~\ref{f:5ionsVCont} (a), continuations on $V$ miss the intermediate region between the red and blue branches. To obtain a complete $I$-$V$ curve, we solve the proposed problem (I2V) with increasing and decreasing $I$. The results are shown in Figure~\ref{f:5ionsVCont} (b), from which we find multiple solutions for some $V$. There are three solutions when $V$ is given in $[V_B, V_E]$ and $[V_D, V_F]$; there are even five solutions when $V$ is given in $[V_E, V_D]$. Here, for instance, $V_B$ denotes the corresponding $V$ value of point $B$. In section \ref{s:TP}, we accurately compute these $V$ values with the method proposed in section \ref{ss:TPMethods}.

\begin{figure}[htbp]
\centering
\includegraphics[scale=.55]{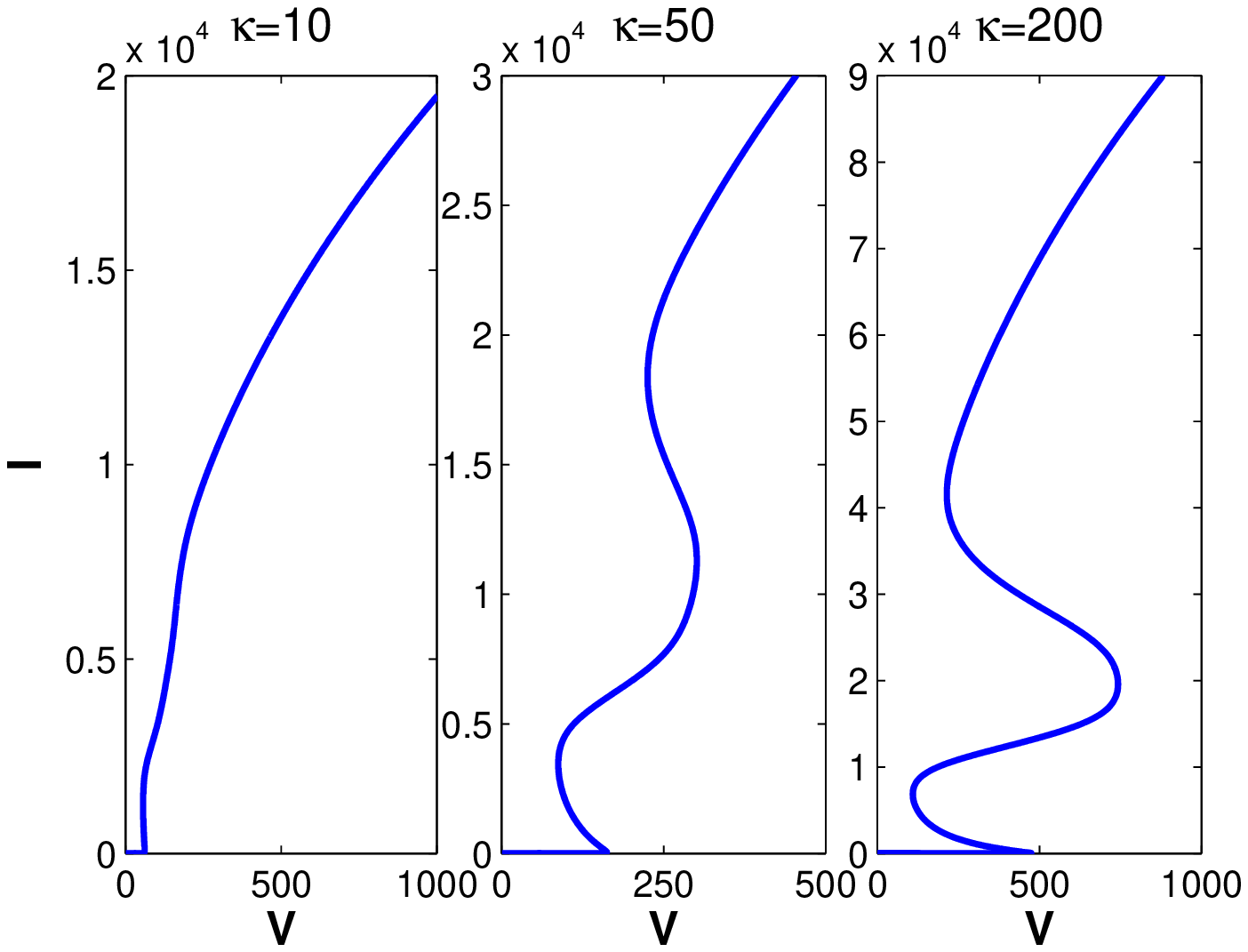}
\includegraphics[scale=.55]{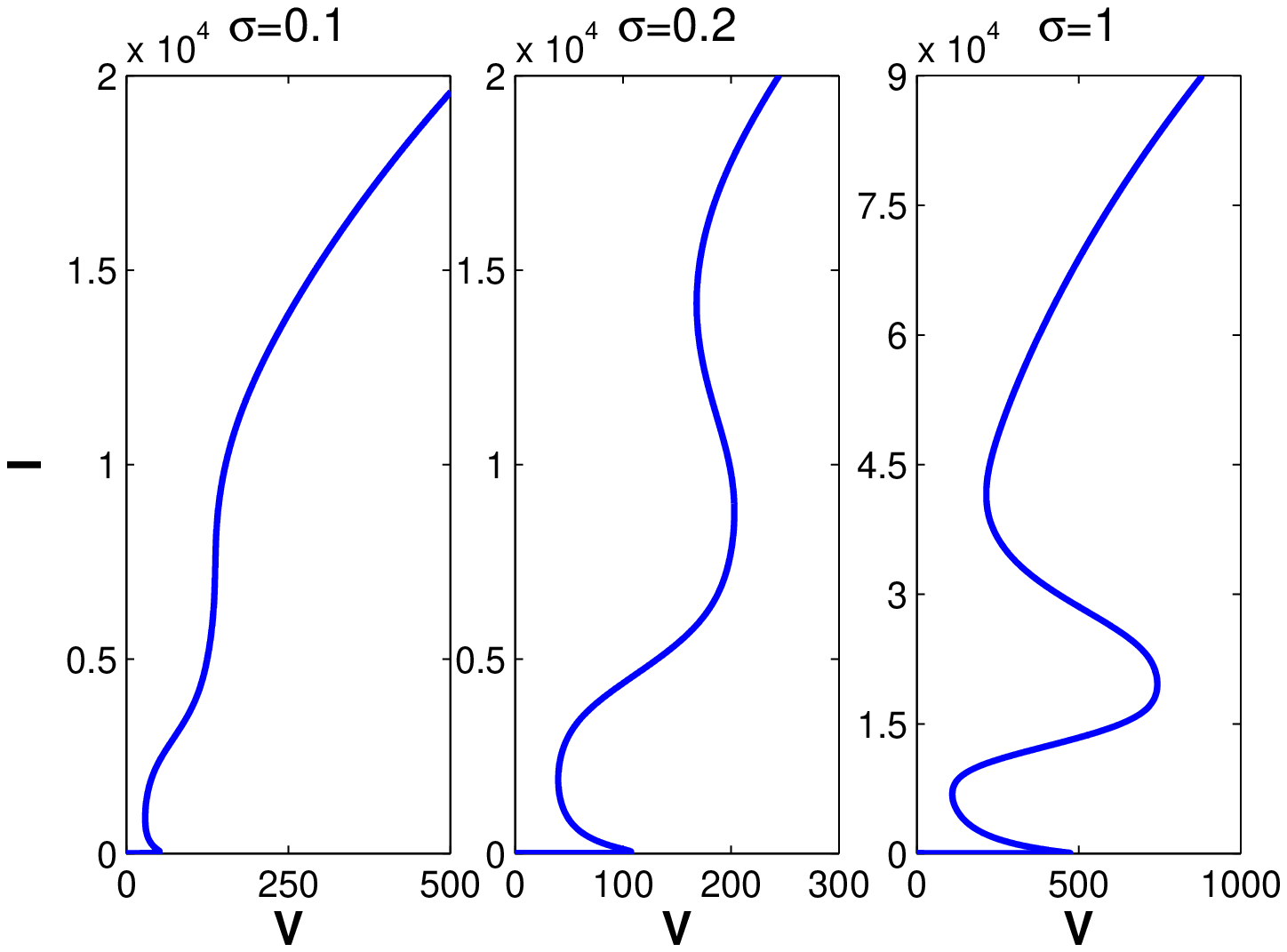}
\caption{ Left: $I$-$V$ curves for different $\kappa$ values with $\sigma=1$; Right:  $I$-$V$ curves for different $\sigma$ values with $\kappa =200$. }
\label{f:DoubleSShape}
\end{figure}
To further study the double S-shaped curve, we are interested in the effect of parameters, $\kappa$ and $\sigma$, on the shape of curves. The results shown in Figure~\ref{f:DoubleSShape} reveal that the shape depends sensitively on $\kappa$ and $\sigma$. We first  set $\sigma=1$ and find that, as $\kappa$ increases from $10$ to $200$, the curve gradually turns from a monotone profile into more and more obvious S shapes. Also, we consider different values of $\sigma$ with a fixed $\kappa=200$. We observe from Figure~\ref{f:DoubleSShape} (b) that the $I$-$V$ curve becomes more and more wavy as $\sigma$ increases, giving rise to multiple (more than three) solutions.

\subsection{Turning Points}\label{s:TP}
\begin{figure}[htbp]
\centering
\includegraphics[scale=.7]{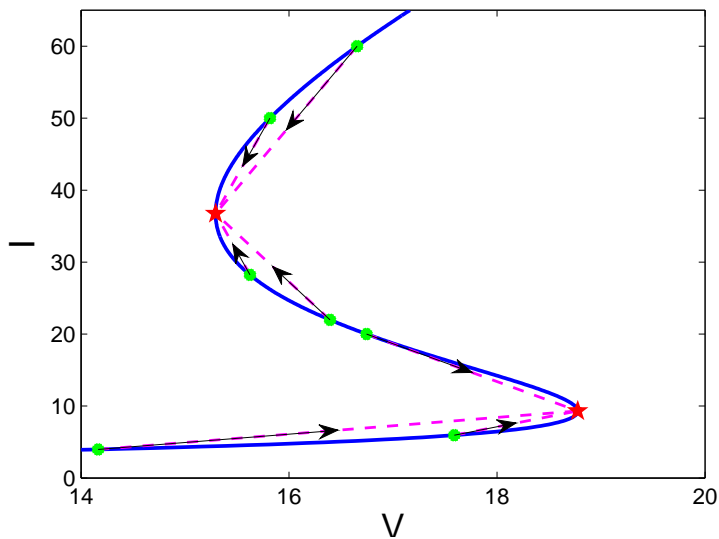}
\caption{Plot for computing turning points on an S-shaped $I$-$V$ curve from different initial guesses, which are labeled by green dots. The red stars represent the turning points found by solving the problem (IVAug).  }
\label{f:2ionsTurningPts}
\end{figure}
In this example, we consider locating turning points on S-shaped and double S-shaped I-V curves. With obtained turning points, where $d V/d I =0$, we are able to determine threshold values over which the system undergoes hysteresis, as well as the interval for $V$ in which the ssPNP equations have multiple solutions.  We remark that there are multiple solutions to the problem (IVAug), corresponding to the multiple turning points on $I$-$V$ curves. To find these solutions, we start from different initial guesses that are obtained by augmenting solutions to the problem (I2V).

Here we report the numerical results from solving the problem (IVAug). In Figure~\ref{f:2ionsTurningPts}, we solve the problem with several initial guesses that are shown by green dots. Such dots are all on the curve, since our initial guesses are obtained based on solutions to the ssPNP equations. Interestingly, these initial guesses, as expected, converge to adjacent turning points. In this case, there are two turning points with coordinates $(15.29, 36.80)$ and $(18.81, 9.27)$ on the $I$-$V$ plane. Therefore, as $V$ increases, the current jumps from the low-current branch to the high-current branch when $V$ becomes larger than $18.81$. Meanwhile, the current jumps from the high-current branch to low-current branch when $V$ gets less than a threshold value $15.29$. Therefore, there are three solutions to the ssPNP equations when the applied voltage belongs to the interval $[15.29, 18.81]$.

We have realized that the threshold value, 15.29 (thermal voltages),  is a bit higher than physiologically relevant values, which often range from $1$ to $6$ (thermal voltages).  However, such a threshold value has sensitive dependence on parameters, especially the fixed charge. The parameters that we used in numerical examples have not been optimized, aiming to lowering the threshold value. It is expected that the threshold value can be lowered to physiologically relevant values by using a certain profile of fixed charges. Optimization methods can be also introduced to lower the threshold value. Similar problems have been studied in the design of bistable optical devices~\cite{Gibbs_Book}. 

\begin{figure}[htbp]
\centering
\includegraphics[scale=.78]{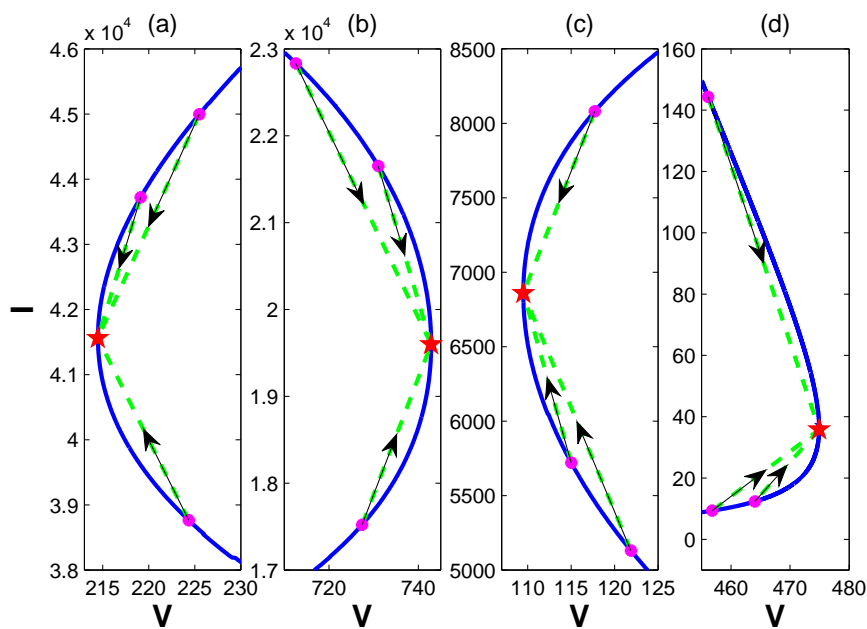}
\caption{Plots for computing four turning points on a double S-shaped $I$-$V$ curve from different initial guesses, which are labeled by pink dots. The red stars represent the turning points found by solving the problem (IVAug).  The double S-shaped  $I$-$V$ curve is shown in the last plot of Figure~\ref{f:5ionsVCont}.}
\label{f:5ionsTurningPts}
\end{figure}

Also, we compute turning points on the double S-shaped $I$-$V$ curve (Figure~\ref{f:5ionsVCont} (b)) discussed in section \ref{s:SSIV}. Although there are $23$ unknowns for the problem (IVAug) with five species of ions, the numerical iterations converge efficiently and robustly to four different turning points. This indicates that our numerical methods are very effective in dealing with multiple species of ions.  As displayed in Figure~\ref{f:5ionsTurningPts}, different initial guesses converge to turning points ssPNP are close. For instance, in Figure~\ref{f:5ionsTurningPts} (a), three different initial guesses converge to the nearest turning point.  There are four turning points on the curve with coordinates $(214.48, 4.16\times10^4)$, $(742.87, 1.96\times10^4)$, $(109.50, 6.86\times10^3)$, and $(474.93, 35.83)$ on the $I$-$V$ plane. As such, there are three solutions in intervals $[109.50, 214.48]$ and $[474.93, 742.87]$, and five solutions in the interval $[214.48, 474.93]$.  Similar to the S-shaped curve, the current has abrupt jumps when the applied voltages are given at $109.50$, $214.48$, $474.93$, and $742.87$, switching the conductance states.

\section{Conclusions}
In this work, we have studied multiple solutions to the steady-state Poisson--Nernst--Planck (ssPNP) equations. We propose four different formulations of the problem to compute complete current-voltage ($I$-$V$) curves and current-concentration relations. In addition, we develop numerical continuations to provide good initial guesses for Newton's iterations that solve nonlinear algebraic equations resulting from discretization. We also have developed a computational method to locate turning points on $I$-$V$ curves, based on an observation that $V$ is a function of $I$ and turing points are achieved at $d V/ d I =0$.  Our numerical results demonstrate that the developed computational methods are robust and effective in solving the  ssPNP equations with multiple ionic species. For instance, the iterations converge robustly and efficiently to turning points on $I$-$V$ curves for channels with five ionic species, in which the corresponding problem has $23$ unknowns.

Of much interest is that we have found voltage- and concentration-induced hysteretic response of current to varying applied voltages and boundary concentrations, respectively.  The results have revealed that conductance states of ionic channels depend on values of applied voltages and boundary concentrations, as well as their history values, showing memory effects. We also study the effect of parameters, e.g., $\kappa$, $\sigma$, and boundary concentrations, on the behavior of such hysteresis.  The abrupt switching of conductance states in hysteresis loops is reminiscent of gating phenomenon, although our model is purely deterministic. Our study may shed some light on the mechanism of gating of ion channels from the viewpoint of deterministic ionic transport.

We have studied multiple solutions to the ssPNP equations with five ionic species. Interestingly, we find a double S-shaped $I$-$V$ curve which gives five solutions to the ssPNP equations, when the applied voltage is given in some interval. As far as we know, the existence of five solutions has been not reported in previous publications. To further explore multiple solutions and hysteresis, we accurately compute turning points on $I$-$V$ curves using a proposed numerical method. With the obtained location of turning points, we are able to determine threshold values over which hysteresis occurs, and the interval for $V$, in which the ssPNP equations have multiple solutions.  These results indicate that the numerical approaches developed here may become useful in studying hysteresis phenomenon.

We now discuss several issues and possible further refinements of our work. First, we have not fully understood the effect of parameters in the ssPNP equations, with multiple ionic species, on the existence of multiple solutions. From our numerical results, we understand that the parameters, $\kappa$, $\sigma$, and boundary concentrations, have strong impact on the shape of $I$-$V$ curves and, therefore, the existence of multiple solutions. Our experience from numerical simulations is that, in addition to the parameters mentioned above, the fixed charge profile with piecewise constants of alternating signs is crucial to the presence of multiple solutions. It is interesting to study that whether the discontinuity or sign alternation is a necessary factor behind the phenomenon of hysteresis.  More rigorous mathematical analysis is therefore needed.

Second, the PNP equations used here do not account for ionic size effect which may become important for crowded environments. The existence of multiple solutions to the PNP theory with size effect has been studied with elegant, rigorous mathematical analysis\cite{LinBob_Nonlinearity15,HungMihn_arXiv15,Gavish_arXiv17}. Note that, in these works, the current has been set to be zero to facilitate the analysis. It is more interesting to study the current-voltage relation when multiple solutions exist. In future, we will apply our developed methods to the PNP theory with size effects.

Finally, it is of great importance to study the stability of the multiple solutions found in our work. In a typical hysteresis loop, the solutions on high-current branch and low-current branch are both stable, but that on the intermediate branch are unstable. Such a bistability phenomenon is ubiquitous in optics\cite{Gibbs_Book}, biological systems\cite{Hyst_Math_Book}, etc. However, it is not easy to study the stability theoretically since there are no closed forms for  the nonhomogeneous multiple solutions found here. Linear stability analysis with the help of numerics will be one of our future works.

 \bigskip
\noindent{\bf Acknowledgments.}
J. Ding and S. Zhou acknowledge supports from Soochow University through a start-up Grant (Q410700415), Natural Science Foundation of Jiangsu Province (BK20160302), and National Natural Science Foundation of China (NSFC 11601361 and NSFC 21773165). The authors also want to thank the anonymous reviewers for their helpful comments.

\bibliographystyle{plain}
\bibliography{PNP}

\begin{thebibliography}{10}

\bibitem{BobLiu_SIADS08}
N.~Abaid, R.~Eisenberg, and W.~Liu.
\newblock Asymptotic expansions of {I}-{V} relations via a
  {P}oisson--{N}ernst--{P}lanck system.
\newblock {\em Phys. Rev. E}, 7:1507--1526, 2008.

\bibitem{Altomare_JGenPhy01}
C.~Altomare, A.~Bucchi, E.~Camatini, M.~Baruscotti, C.~Viscomi, A.~Moroni, and
  D.~DiFrancesco.
\newblock Integrated allosteric model of voltage gating of {HCN} channels.
\newblock {\em J. Gen. Physiol.}, 117:519--532, 2001.

\bibitem{Andersson_MathBiosci2010}
T.~Andersson.
\newblock Exploring voltage-dependent ion channels in silico by hysteretic
  conductance.
\newblock {\em Math. Biosci.}, 226:16--27, 2010.

\bibitem{BarChenBob_SIAP92}
V.~Barcilon, D.~Chen, and R.~Eisenberg.
\newblock Ion flow through narrow membrane channels: Part {II}.
\newblock {\em SIAM J. Appl. Math.}, 52:1405--1425, 1992.

\bibitem{BarChenBob_SIAP97}
V.~Barcilon, D.~Chen, R.~Eisenberg, and J.~Jerome.
\newblock Qualitative properties of steady-state {P}oisson--{N}ernst--{P}lanck
  systems: {P}erturbation and simulation study.
\newblock {\em SIAM J. Appl. Math.}, 57:631--648, 1997.

\bibitem{BazantChuBayly_SIAP06}
M.~Bazant, K.~Chu, and B.~Bayly.
\newblock Current-voltage relations for electrochemical thin films.
\newblock {\em SIAM J. Appl. Math.}, 65:1463--1484, 2006.

\bibitem{Hyst_Math_Book}
M.~Brokate and J.~Spreckels.
\newblock {\em Hysteresis and Phase Transition, in: Applied Mathematical
  Sciences}.
\newblock Springer, New York, 1996.

\bibitem{Hyst_RevModPhys99}
B.~Chakrabarti and M.~Acharyya.
\newblock Dynamic transitions and hysteresis.
\newblock {\em Rev. Mod. Phys.}, 71:1--30, 1999.

\bibitem{Das_PRE2012}
B.~Das, K.~Banerjee, and G.~Gangopadhyay.
\newblock Entropy hysteresis and nonequilibrium thermodynamic efficiency of ion
  conduction in a voltage-gated potassium ion channel.
\newblock {\em Phys. Rev. E}, 86:061915, 2013.

\bibitem{BobHyonLiu_JCP10}
B.~Eisenberg, Y.~Hyon, and C.~Liu.
\newblock Energy variational analysis {EnVarA} of ions in water and channels:
  Field theory for primitive models of complex ionic fluids.
\newblock {\em J. Chem. Phys.}, 133:104104, 2010.

\bibitem{BergLiu_SIMA07}
B.~Eisenberg and W.~Liu.
\newblock {P}oisson--{N}ernst--{P}lanck systems for ion channels with permanent
  charges.
\newblock {\em SIAM J. Math. Anal.}, 38:1932--1966, 2007.

\bibitem{BobLiuXu_Nonlinearity15}
B.~Eisenberg, W.~Liu, and H.~Xu.
\newblock Reversal permanent charge and reversal potential: case studies via
  classical {P}oisson--{N}ernst--{P}lanck models.
\newblock {\em Nonlinearity}, 28:103--127, 2015.

\bibitem{Noble_PTRSocA_09}
M.~Fink and D.~Noble.
\newblock Markov models for ion channels: {V}ersatility versus identifiability
  and speed.
\newblock {\em Phil. Trans. R. Soc. A}, 367:2161--2179, 2009.

\bibitem{FologeaBBActa_2011}
D.~Fologea, E.~Krueger, Y.~Mazur, C.~Stith, Y.~Okuyama, R.~Henry, and
  G.~Salamo.
\newblock Bi-stability, hysteresis, and memory of voltage-gated lysenin
  channels.
\newblock {\em Biochim. Biophys. Acta}, 1808:2933--2939, 2011.

\bibitem{Gavish_arXiv17}
N.~Gavish.
\newblock {P}oisson--{N}ernst--{P}lanck equations with steric effects --
  non-convexity and multiple stationary solutions.
\newblock {\em arXiv preprint}, 1:07164, 2017.

\bibitem{Gibbs_Book}
H.~Gibbs.
\newblock {\em Optical Bistability: Controlling Light with Light}.
\newblock Academic, New York, 1985.

\bibitem{Usherwood_Nat81}
K.~Gration, J.~Lambert, R.~Ramsey, and P.~Userwood.
\newblock Non-random openings and concentration-dependent lifetimes of
  {G}lutamate-gated channels in muscle membrane.
\newblock {\em Nature}, 291:423--425, 1981.

\bibitem{Hille_Book2001}
B.~Hille.
\newblock {\em Ion Channels of Excitable Membranes}.
\newblock Sinauer Associates, 3rd edition, 2001.

\bibitem{Hodgkin_PRSL58}
A.~Hodgkin.
\newblock Ionic movements and electrical activity in giant nerve fibres.
\newblock {\em Proc. R. Soc. Lond. B Biol. Sci.}, 148:1--37, 1958.

\bibitem{HH2_JPhys52}
A.~Hodgkin and A.~Huxley.
\newblock The components of membrane conductance in the giant axon of loligo.
\newblock {\em J. Physiol.}, 116:473--496, 1952.

\bibitem{HH1_JPhys52}
A.~Hodgkin and A.~Huxley.
\newblock Currents carried by sodium and potassium ions through the membrane of
  the giant axon of loligo.
\newblock {\em J. Physiol.}, 116:449--472, 1952.

\bibitem{HH3_JPhys52}
A.~Hodgkin and A.~Huxley.
\newblock A quantitative description of membrane current and its application to
  conduction and excitation in nerve.
\newblock {\em J. Physiol.}, 117:500--544, 1952.

\bibitem{HungMihn_arXiv15}
L.~Hung and M.~Minh.
\newblock Stationary solutions to the {P}oisson--{N}ernst--{P}lanck equations
  with steric effects.
\newblock {\em arXiv preprint}, 1:02456, 2015.

\bibitem{HyonLiuBob_CMS10}
Y.~Hyon, B.~Eisenberg, and C.~Liu.
\newblock A mathematical model for the hard sphere repulsion in ionic
  solutions.
\newblock {\em Commun. Math. Sci.}, 9:459--475, 2010.

\bibitem{HyonLiuBob_JPCB12}
Y.~Hyon, C.~Liu, and B.~Eisenberg.
\newblock {PNP} equations with steric effects: a model of ion flow through
  channels.
\newblock {\em J. Phys. Chem. B}, 116:11422--11441, 2012.

\bibitem{JiLiuZhang_SIAP15}
S.~Ji, W.~Liu, and M.~Zhang.
\newblock Effects of (small) permanent charges and channel geometry on ionic
  flows via classical {P}oisson--{N}ernst--{P}lanck models.
\newblock {\em SIAM J. Appl. Math.}, 75:114--135, 2015.

\bibitem{JiaLIuZhang_DCDSB16}
Y.~Jia, W.~Liu, and M.~Zhang.
\newblock Qualitative properties of ionic flows via
  {P}oisson--{N}ernst--{P}lanck systems with bikerman's local hard-sphere
  potential: Ion size effects.
\newblock {\em Discrete Contin. Dyn. Syst.}, 21:1775--1802, 2016.

\bibitem{Cui_PflugersArch94}
I.~Josephson and Y.~Cui.
\newblock Voltage- and concentration-dependent effects of lidocaine on cardiac
  {N}a channel and {C}a channel gating charge movements.
\newblock {\em Pflugers Arch.}, 428:485--491, 1994.

\bibitem{BVP4C}
J.~Kierzenka and L.~Shampine.
\newblock A {BVP} solver based on residual control and the {M}atlab {PSE}.
\newblock {\em ACM Trans. Math. Softw.}, 27:299--316, 2001.

\bibitem{BazantSteric_PRE07}
M.~Kilic, M.~Bazant, and A.~Ajdari.
\newblock Steric effects in the dynamics of electrolytes at large applied
  voltages. {II}. {M}odified {P}oisson--{N}ernst--{P}lanck equations.
\newblock {\em Phys. Rev. E}, 75:021503, 2007.

\bibitem{Krueger_BiophyChem13}
E.~Krueger, R.~Faouri, D.~Fologea, R.~Henry, D.~Straub, and G.~Salamo.
\newblock A model for the hysteresis observed in gating of lysenin channels.
\newblock {\em Biophys. Chem.}, 184:126--130, 2013.

\bibitem{LeeHyonLinLiu_Nonlinearity11}
C.~Lee, H.~Lee, Y.~Hyon, T.~Lin, and C.~Liu.
\newblock New poisson-boltzmann type equations: one-dimensional solutions.
\newblock {\em Nonlinearity}, 24:431--458, 2011.

\bibitem{LiLiuXuZhou_Nonliearity13}
B.~Li, P.~Liu, Z.~Xu, and S.~Zhou.
\newblock Ionic size effects: generalized boltzmann distributions, counterion
  stratification, and modified debye length.
\newblock {\em Nonlinearity}, 26(10):2899, 2013.

\bibitem{LiWenZhou_CMS16}
B.~Li, J.~Wen, and S.~Zhou.
\newblock Mean-field theory and computation of electrostatics with ionic
  concentration dependent dielectrics.
\newblock {\em Commun. Math. Sci.}, 14:249--271, 2016.

\bibitem{LinLiuZhang_SIADS13}
G.~Lin, W.~Liu, Y.~Yi, and M.~Zhang.
\newblock {P}oisson--{N}ernst--{P}lanck systems for ion flow with density
  functional theory for local hard-sphere potential.
\newblock {\em SIAM J. Appl. Dyn. Syst.}, 12:1613--1648, 2013.

\bibitem{LinBob_CMS14}
T.~Lin and B.~Eisenberg.
\newblock A new approach to the {Lennard-Jones} potential and a new model:
  {PNP}-steric equations.
\newblock {\em Commun. Math. Sci.}, 12:149--173, 2014.

\bibitem{LinBob_Nonlinearity15}
T.~Lin and B.~Eisenberg.
\newblock Multiple solutions of steady-state {P}oisson--{N}ernst--{P}lanck
  equations with steric effects.
\newblock {\em Nonlinearity}, 28:2053--2080, 2015.

\bibitem{WLiu_SIAP05}
W.~Liu.
\newblock Geometric singular perturbation approach to steady-state
  {P}oisson--{N}ernst--{P}lanck systems.
\newblock {\em SIAM J. Appl. Math.}, 65:754--766, 2005.

\bibitem{WLiu_JDE09}
W.~Liu.
\newblock One-dimensional steady-state {P}oisson--{N}ernst--{P}lanck systems
  for ion channels with multiple ion species.
\newblock {\em J. Differ. Equ.}, 246:428--451, 2009.

\bibitem{LiuXu_JDE15}
W.~Liu and H.~Xu.
\newblock A complete analysis of a classical {P}oisson--{N}ernst--{P}lanck
  model for ionic flow.
\newblock {\em J. Differ. Equ.}, 258:1192--1228, 2015.

\bibitem{BZLu_BiophyJ11}
B.~Lu and Y.~Zhou.
\newblock {Poisson-Nernst-Planck} equations for simulating biomolecular
  diffusion-reaction processes {II}: Size effects on ionic distributions and
  diffusion-reaction rates.
\newblock {\em Biophys. J.}, 100:2475--2485, 2011.

\bibitem{MaXu_JCP14}
M.~Ma and Z.~Xu.
\newblock Self-consistent field model for strong electrostatic correlations and
  inhomogeneous dielectric media.
\newblock {\em J. Chem. Phys.}, 141:244903, 2014.

\bibitem{MacKinnon04_ACIE04}
R.~MacKinnon.
\newblock Potassium channels and the atomic basis of selective ion conduction.
\newblock {\em Angew. Chem. Int. Ed. Engl.}, 43:4265--4277, 2004.

\bibitem{Roope_JGPhys05}
R.~Mannikko, S.~Pandey, H.~Larsson, and F.~Elinder.
\newblock Hysteresis in the voltage dependence of {HCN} channels: Conversion
  between two modes affects pacemaker properties.
\newblock {\em J. Gen. Physiol.}, 125:305--326, 2005.

\bibitem{PMarkowich_Book}
P.~Markowich.
\newblock {\em The Stationary Semiconductor Device Equations}.
\newblock Springer-Verlag, New York, 1986.

\bibitem{MockExample}
M.~Mock.
\newblock An example of nonuniqueness of stationary solutions in semiconductor
  device models.
\newblock {\em COMPEL}, 1:165--174, 1982.

\bibitem{Nache_NatComm2013}
V.~Nache, T.~Eick, E.~Schulz, R.~Schmauder, and K.~Benndorf.
\newblock Hysteresis of ligand binding in {CNGA2} ion channels.
\newblock {\em Nat. Commun.}, 4:1--9, 2013.

\bibitem{Bezrukov_JCP06}
M.~Pustovoit, A.~Berezhkovskii, and S.~Bezrukova.
\newblock Analytical theory of hysteresis in ion channels: Two-state model.
\newblock {\em J. Chem. Phys.}, 125, 2006.

\bibitem{BZLu_JSP16}
Y.~Qiao, X.~Liu, M.~Chen, and B.~Lu.
\newblock A local approximation of fundamental measure theory incorporated into
  three dimensional {Poisson--Nernst--Planck} equations to account for hard
  sphere repulsion among ions.
\newblock {\em J. Stat. Phys.}, 163:156--174, 2016.

\bibitem{BZLu_JCP14}
Y.~Qiao, B.~Tu, and B.~Lu.
\newblock Ionic size effects to molecular solvation energy and to ion current
  across a channel resulted from the nonuniform size-modified {PNP} equations.
\newblock {\em J. Chem. Phys.}, 140:174102, 2014.

\bibitem{Bezrukov_EBioPhyJ15}
S.~Rappaport, O.~Teijido, D.~Hoogerheide, T.~Rostovtseva, A.~Berezhkovskii, and
  S.~Bezrukov.
\newblock Conductance hysteresis in the voltage-dependent anion channel.
\newblock {\em Euro. Biophys. J.}, 44:465--472, 2015.

\bibitem{Yamoah_BioPhyJ03}
A.~Rodriguez-Contreras and E.~Yamoah.
\newblock Effects of permeant ion concentrations on the gating of {L}-type
  {$Ca^{2+}$} channels in hair cells.
\newblock {\em Biophys. J.}, 84:3457--3459, 2003.

\bibitem{Rubinstein_SIAP87}
I.~Rubinstein.
\newblock Multiple steady state solutions in one-dimensional electrodiffusion
  with local electroneutrality.
\newblock {\em SIAM J. Appl. Math.}, 47:1076--1093, 1987.

\bibitem{RubinsteinSIAM_Book}
I.~Rubinstein.
\newblock {\em Electro-Diffusion of Ions}.
\newblock SIAM Stud. Appl. Math., Philadelphia, 1990.

\bibitem{SchussNadlerBob_PRE01}
Z.~Schuss, B.~Nadler, and R.~Eisenberg.
\newblock Derivation of {Poisson} and {Nernst}--{Planck} equations in a bath
  and channel from a molecular model.
\newblock {\em Phys. Rev. E}, 64:1--14, 2001.

\bibitem{NeuroBook}
G.~Shepherd.
\newblock {\em Neurobiology}.
\newblock Oxford University Press, Oxford, 1994.

\bibitem{Sigworth_QRB94}
F.~Sigworth.
\newblock Voltage gating of ion channels.
\newblock {\em Q. Rev. Biophys.}, 27:1--40, 1994.

\bibitem{SingerGillBob_ESIAM08}
A.~Singer, J.~Norbury D.~Gillespie, and R.~Eisenberg.
\newblock Singular perturbation analysis of the steady-state
  {P}oisson--{N}ernst--{P}lanck system: Applications to ion channels.
\newblock {\em Euro. J. Appl. Math.}, 19:541--560, 2008.

\bibitem{Steinruck_SIAP89}
H.~Steinruck.
\newblock A bifurcation analysis of the one-dimensional steady state
  semiconductor device equations.
\newblock {\em SIAM J. Appl. Math.}, 49:1101--1121, 1989.

\bibitem{Cuello_PNAS17}
C.~Tilegenovaa, D.~Cortesa, and L.~Cuello.
\newblock Hysteresis of {KcsA} potassium channel's activation-deactivation
  gating is caused by structural changes at the channel's selectivity filter.
\newblock {\em Proc. Natl. Acad. Sci. USA.}, 114:3234--3239, 2017.

\bibitem{WangHeHuang_PRE14}
X.~Wang, D.~He, J.~Wylie, and H.~Huang.
\newblock Singular perturbation solutions of steady-state
  {P}oisson--{N}ernst--{P}lanck systems.
\newblock {\em Euro. J. Appl. Math.}, 89:022722, 2014.

\bibitem{Ward_SIAP91}
M.~Ward, L.~Reyna, and F.~Odeh.
\newblock Multiple steady-state solutions in a multijunction semiconductor
  device.
\newblock {\em SIAM J. Appl. Math.}, 51:90--123, 1991.

\bibitem{XuShengLiu_CMS14}
S.~Xu, P.~Sheng, and C.~Liu.
\newblock An energetic variational approach for ion transport.
\newblock {\em Commun. Math. Sci.}, 12:779--789, 2014.

\bibitem{XuMaLiu_PRE14}
Z.~Xu, M.~Ma, and P.~Liu.
\newblock Self-energy-modified {P}oisson--{N}ernst--{P}lanck equations: {WKB}
  approximation and finite-difference approaches.
\newblock {\em Phys. Rev. E}, 90:013307, 2014.

\bibitem{IonChanel_HandbookCRC15}
J.~Zheng and M.~Trudeau.
\newblock {\em Handbook of ion channels}.
\newblock CRC Press, 2015.

\bibitem{ZhouWangLi_PRE11}
S.~Zhou, Z.~Wang, and B.~Li.
\newblock Mean-field description of ionic size effects with non-uniform ionic
  sizes: {A} numerical approach.
\newblock {\em Phys. Rev. E}, 84:021901, 2011.

\end{thebibliography}
\end{document}